\magnification 1200
\input amssym.def
\input amssym.tex
\parindent = 40 pt
\parskip = 12 pt
\font \heading = cmbx10 at 14 true pt
 at 22 true pt
\font \medheading =cmbx7 at 16 true pt
\font \small = cmr7 at 10 true pt
\def \R{{\bf R}}

\centerline{\medheading $L^p$ Sobolev regularity of averaging operators over }
\centerline{\medheading   hypersurfaces and the Newton polyhedron}
\rm
\line{}
\line{}
\centerline{\heading Michael Greenblatt}
\line{}
\centerline{February 18, 2018}
\line{}
\line{}\baselineskip = 4 pt
 \line{}\line{}\line{}\line{}
 {\narrower \noindent \small $L^p$ to $L^p_{\beta}$ boundedness results are proven for translation invariant averaging operators
 over hypersurfaces in Euclidean space. The operators can either be Radon transforms or averaging operators with multiparameter fractional integral kernel. In many cases, the amount $\beta > 0 $ of smoothing
proven is optimal up to endpoints, and in such situations this amount of smoothing can be computed explicitly through the use of appropriate Newton polyhedra. \par}
\baselineskip = 12 pt
\font \heading = cmbx10 at 14 true pt
\line{}
\line{}
\line{}
\noindent{\heading 1. Introduction and theorem statements.} 

In this paper we consider convolution operators with hypersurface measures on $\R^{n+1}$
of the following form, where $\bf{x}$ denotes $(x_1,...,x_n)$ and $\bf{t}$ denotes $(t_1,...,t_n)$.
$$Tf({\bf x}, x_{n+1}) = \int_{\R^n} f({\bf x} - {\bf t}, x_{n+1} - S({\bf t})) K({\bf t})\,d{\bf t}\eqno (1.1)$$
Here $S({\bf t})$ is a real-analytic function on a neighborhood $U$ of the origin and $K({\bf t})$ is a function supported in $U$ that
is $C^1$ on $\{{\bf t} \in U: t_i \neq 0$ for all $i\}$ which satisfies estimates as follows. Write ${\bf t} = ({\bf t}_1,...,{\bf t}_m)$,
where ${\bf t}_i$ denotes  $(t_{i1},...,t_{il_i})$ such that the various $t_{ij}$ variables comprise the whole list $t_1,...,t_n$. Then for some
$0 \leq \alpha_i < l_i$ and some $C > 0$ we assume the following estimates.
$$|K({\bf t})| \leq C\prod_{k=1}^m |{\bf t}_k|^{-\alpha_k} \eqno (1.2a) $$
$$|\partial_{t_{ij}}K({\bf t})| \leq C {1 \over |t_{ij}|} \prod_{k=1}^m |{\bf t}_k|^{-\alpha_k} {\hskip 0.4 in} {\rm \,\,for\,\,all\,\,} i {\rm\,\, and\,\,} j \eqno (1.2b) $$
In the parlance of harmonic analysis, operators satisfying $(1.1), (1.2a), (1.2b)$ are sometimes referred to as fractional Radon transforms or fractional singular 
Radon transforms. The case where each $\alpha_i = 0$ includes the traditional Radon transform operators.
Our goal will be to prove $L^p(\R^{n+1})$ to $L^p_{\beta}(\R^{n+1})$ boundedness properties for $T$ that are sharp up 
to endpoints and which can be computed explicitly with the help of certain Newton polyhedra, as will be described after the proof of Lemma 2.1. The author's earlier 
paper [G1] covers, among other things, the situation where $p = 2$ and each $\alpha_i = 0$, and our methods and results will extend those of [G1].

By the translation and rotation invariance properties 
of convolution operators, without loss of generality we may assume that
$$S(0,...,0) = 0 {\hskip 0.65 in} \nabla S(0,...,0) = (0,...,0) \eqno (1.3)$$
To avoid trivialities, we also assume $S$ is not identically zero. In order to describe our results, we will make use of the following terminology.

\noindent {\bf Definition 1.1.} Let $f({\bf t})$ be a real analytic function defined on a neighborhood of the origin in 
$\R^n$, and
let $f({\bf t}) = \sum_{\alpha} f_{\alpha}{\bf t}^{\alpha}$ denote the Taylor expansion of $f({\bf t})$ at the origin.
For any $\alpha$ for which $f_{\alpha} \neq 0$, let $Q_{\alpha}$ be the octant $\{{\bf t} \in \R^n: 
t_i \geq \alpha_i$ for all $i\}$. Then the {\it Newton polyhedron} $N(f)$ of $f({\bf t})$ is defined to be 
the convex hull of all $Q_{\alpha}$.  

A Newton polyhedron can contain faces of various dimensions in various configurations. 
These faces can be either compact or unbounded. A vertex of $N(f)$ is considered to be a face of dimension zero.

\noindent {\bf Definition 1.2.} Where $f({\bf t})$ is as in Definition 1.1, we define $f^*({\bf t})$ by 
$$f^*({\bf t}) = \sum_{(v_1,...,v_n)\,\,a \,\,vertex \,\,of\,\,N(f)} |t_1^{v_1}...\,t_n^{v_n}| \eqno (1.4)$$
By Lemma 2.1 of [G2], there
is a constant $C$ such that for all ${\bf t}$ one has $|f({\bf t})| \leq Cf^*({\bf t})$. 

Let $d\mu$ denote the measure $\prod_{k=1}^m |{\bf t}_k|^{-\alpha_k}\,dm$, where $m$ denotes Lebesgue measure.
In Lemma 2.1, we will show that there is an $r_0 > 0$, an $a_0> 0$, and an integer $d_0$ satisfying $0 \leq d_0 \leq n-1$, such that if $r < r_0$ then there are positive constants $b_r$ and $B_r$ such that for $0 < \epsilon < {1 \over 2}$ we have
$$b_r \epsilon^{a_0} |\ln \epsilon|^{d_0} < \mu(\{{\bf t} \in (0,r)^n:  S^*({\bf t}) < \epsilon\}) < B_r \epsilon^{a_0} |\ln\epsilon|^{d_0}\eqno (1.5)$$
The quantity $a_0$ will play a key role in our results. In order to state our main theorem, we will also need the following definitions.

\noindent {\bf Definition 1.3.} Suppose $F$ is a compact face of the Newton polyhedron $N(f)$. Then
if $f({\bf t}) = \sum_{\alpha} f_{\alpha}{\bf t}^{\alpha}$ denotes the Taylor expansion of $f$ like above, 
define $f_F({\bf t}) = \sum_{\alpha \in F} f_{\alpha}{\bf t}^{\alpha}$.

\noindent {\bf Definition 1.4.} For $f({\bf t})$ as above, we denote by $o(f)$ the maximum order of any zero of any $f_F({\bf t})$ on $(\R - \{0\})^n$. We take
$o(f) = 0$ if there are no such zeroes.

\noindent We now come to our main theorem. 

 \noindent {\bf Theorem 1.1.} Suppose $S({\bf t})$ is a real analytic function on a neighborhood of the origin satisfying $(1.3)$. Let $g = 
 \min(a_0, l_1 - \alpha_1,...,\l_m - \alpha_m)$, where the $\alpha_i$ and $l_i$ are as in the beginning of the paper and $a_0$ is as in $(1.5)$. Then there is a neighborhood
 $V$ of the origin such that if $K({\bf t})$ is supported on $V$ and satisfies $(1.2a)-(1.2b)$ then the following hold.
 
 \noindent {\bf 1)} Let $A$ denote the open triangle with vertices $({1 \over 2}, {1 \over \max(o(S), 2)})$, $(0,0)$, and $(1,0)$, and  let $B = \{(x,y) \in A: 
 y < g\}$. Then $T$ is bounded from $L^p(\R^{n+1})$ to $L^p_{\beta}(\R^{n+1})$ if $({1 \over p}, \beta) \in B$. 
 
 \noindent {\bf 2)} Suppose $g < 1$, $K({\bf t})$ is nonnegative, and there exists a positive constant $C_0$ and a neighborhood $N_0$ of the origin such that $K({\bf t}) > C_0\prod_{k=1}^m |{\bf t}_k|^{-\alpha_k}$ on $\{{\bf t} \in N_0: t_i \neq 0$ for all $i\}$. Then if $1 < p < \infty$ and 
 $T$ is bounded from $L^p(\R^{n+1})$ to $L^p_{\beta}(\R^{n+1})$ we must have $\beta \leq g$.

Note that when $g < {1 \over \max(o(S), 2)}$, then the two parts of Theorem 1.1 combined say that for ${1 \over p} \in ( {\max(o(S), 2) \over 2}g, 1 -  {\max(o(S), 2) \over 2}g)$, the amount of $L^p$ 
Sobolev smoothing given by part 1, $g$ derivatives, is optimal except possibly missing the endpoint $\beta = g$. When $g = {1 \over \max(o(S), 2)}$ the same is true for $p = 2$. A natural question to ask is when does the endpoint $\beta = g$ also hold. Although we will not show it here, it
 turns out that when $p = 2$ sometimes this endpoint estimate holds and sometimes it does not. The author does not know what happens in the 
 $p \neq 2$ situation.
 
 It is not hard to show that part 2) of Theorem 1.1 does not necessarily hold when $g \geq 1$. A simple counterexample is the case where $K({\bf t})$ is smooth and nonnegative with $K(0) > 0$,  each 
 $l_i = 1$, each $\alpha_i = 0$,  and  $S({\bf t}) = \sum_{i=1}^n t_i^2$. Then $g = 1$, but by the well-known Fourier transform decay estimates for
 nondegenerate hypersurface measures (see p. 348 of [S]) one has that $T$ is bounded from $L^2(\R^{n+1})$ to $L^2_{{n \over 2}}(\R^{n+1})$.

\noindent {\bf Example 1.} Suppose we are in the case where each $\alpha_i = 0$, such as in the case of smooth $K({\bf t})$. Then $g = 
\min(a_0, l_1,...,l_m)$. Because the triangle $A$ of Theorem 1.1 has its upper vertex no higher than $({1 \over 2}, {1 \over 2})$, if $a_0 > {1 \over 2}$ 
the set $B$ of Theorem 1.1 is just the triangle $A$ and the theorem does not necessarily give any sharp estimates, while if $a_0 \leq  {1 \over 2}$ we have $g = a_0$ and Theorem 1.1 gives a result that is sharp
up to endpoints if $a_0 \leq {1 \over \max(o(S), 2)}$. The $p = 2$ case of these facts follow from Theorem 1.5 of [G1], which in turn generalizes 
portions of a paper of Varchenko [V] that covers the cases where $o(S) = 0$ or $1$. 

As will be explained in greater detail in the proof of Lemma 2.1,
in this example $a_0$ is the reciprocal of  the Newton distance of $S^{*}({\bf t})$, defined as follows.

\noindent {\bf Definition 1.3.} The {\it Newton distance} $d(f)$ of  a function $f({\bf t})$ is defined to be $\inf \{c: c,c,...,c,c) \in N(f)\}$.

\noindent {\bf Example 2.} Suppose $m = 1$, so that $d\mu = |{\bf t}|^{-\alpha_1}\,dm$ for some $\alpha_1 < n$. Then given that $(1.3)$ holds, 
one has $S^*({\bf t}) \leq \sum_{i=1}^n |t_i|^2 = |{\bf t}|^2$, and as a result if $0 < r < 1$ and $\epsilon^{1 \over 2} < r$ one has
$$\mu(\{{\bf t} \in (0,r)^n:  S^*({\bf t}) < \epsilon\}) \geq \mu(\{{\bf t} \in (0,r)^n: |{\bf t}|^2 < \epsilon\})$$
$$= \int_{|{\bf t}| < \epsilon^{1 \over 2}}| {\bf t}|^{-\alpha_1} \,d{\bf t}\eqno (1.6)$$
The integral $(1.6)$ is easily seen to be of the form $C\epsilon^{n - \alpha_1 \over 2}$. Hence $a_0 \leq {n - \alpha_1 \over 2}$. Since $g = \min(a_0,
n - \alpha_1)$, we have that $g = a_0$  whenever $m = 1$. The value of $a_0$ can then be computed using Newton polyhedra as in the proof of Lemma 2.1.

\noindent {\bf Example 3.} Suppose $m = n$ so that each ${\bf t}_k$ is one dimensional. Then we have
$$\mu(\{{\bf t} \in (0,r)^n:  S^*({\bf t}) < \epsilon\}) = \int_{\{{\bf t}\in (0,r)^n:\,  S^*({\bf t}) < \epsilon\}}\prod_{k=1}^n t_k^{-\alpha_k}\,d{\bf t}  \eqno (1.7)$$
We change variables from $t_k$ to $u_k = t_k^{1 - \alpha_k}$. If $S^*({\bf t}) = \sum_{i = 1}^N |t_1|^{v_{i1}}...|t_n|^{v_{in}}$, then we denote by $S^{**}({\bf u})$ the function $\sum_{i = 1}^N |u_1|^{v_{i1} \over 1 - \alpha_1}...
|u_n|^{v_{in} \over 1 - \alpha_n}$. Then for an appropriate box $B$, $(1.7)$ becomes of the form
$$\mu(\{{\bf t} \in (0,r)^n: S^*({\bf t}) < \epsilon\}) = C \int_{\{{\bf u}\in B: \, S^{**}({\bf u}) < \epsilon\}}1 \,d{\bf u}$$
$$= C m(\{{\bf u}\in B:\,  S^{**}({\bf u}) < \epsilon\}) \eqno (1.8)$$
As in Example 1, the reciprocal of the  Newton distance of $S^{**}$ gives the exponent in the growth rate of
 $m(\{{\bf u}\in B:  S^{**}({\bf u}) < \epsilon\})$,  which
by $(1.8)$ then gives the value of $a_0$. 

\noindent {\bf Some history.} 

There has been a great deal of work done on function space boundedness properties of Radon transforms and related operators, so we focus on the
$L^p$ to $L^p_{\beta}$ boundedness questions being considered here. The case of translation-invariant Radon transforms with smooth density functions over curves in $\R^2$ 
was thoroughly analyzed in [Gr] and [C]. In the general non-translation-invariant case for curves in $\R^2$ , again with smooth density functions,  comprehensive $L^p_{\alpha}$ to $L^q_{\beta}$ estimates that are sharp up to endpoints are proven in [Se].

For translation-invariant Radon transforms over two dimensional surfaces in $\R^3$, there are a number of results. If $p = 2$, the level of Sobolev space
smoothing directly translates into a surface measure Fourier transform decay rate problem, and for the case of smooth density functions the stability 
theorems of Karpushkin [Ka1] [Ka2] combined with [V] give some sharp decay rate results. Generalizations to smooth phases then follow from 
[IKM]. When one has a singular density function, in [G3] the author proved some theorems  corresponding to the case $m = 1$ of this paper
that include results that go beyond what is proved here. Other results for singular density functions appear in [G4].

For higher dimensional hypersurfaces, if the density functions are singular enough in the sense that the $\alpha_i$ are close enough to $l_i$, then there will
be an interval $I$ containing 2 such that sharp $L^p$ Sobolev smoothing estimates follow from the results of [St1] when $p \in I$. This extends earlier
work of the author [G5]. If one lets the 
$\alpha_i$ actually equal $l_i$ and one adds an appropriate cancellation condition one has a multiparameter singular Radon transform, and $L^p$ 
boundedness results for such operators were proven in [St1] [St2], extending the results in [CNSW]. Additional results for higher dimensional
hypersurfaces appear in [Cu], and in [PSe] Sobolev space estimates are proven for translation-invariant Radon transforms over curves.

\noindent {\heading 2. Some useful lemmas.} 

\noindent {\bf Lemma 2.1.}  There is an $r_0 > 0$, an $a_0> 0$, and an integer $d_0$ satisfying $0 \leq d_0 \leq n-1$, such that if $r < r_0$ then there are positive constants $b_r$ and $B_r$ such that for $0 < \epsilon < {1 \over 2}$ equation $(1.5)$ holds.

\noindent {\bf Proof.} Given $r > 0$, we divide $(0,r)^n$ (up to a set of measure zero) into $2^{n!}$ regions $\{A_l\}_{l=1}^{2^{n!}}$, where 
each $A_l$ is a region of the form $\{{\bf t} \in (0,r)^n: t_{j_1} < t_{j_2} < ... < t_{j_n}\}$, where  $t_{j_1},...,t_{j_n}$ is a permutation of
the $t$ variables. We focus our attention one one such $A_l$ and we let $u_i$ denote $t_{j_i}$. Then in the $u_i$ variables the function $S^*({\bf t})$ becomes
a function $S_l^*({\bf u})$ of the form $\sum_{i=1}^N |u^{w_i}|$, where the components of a given $w_i$ are a permutation of that of $v_i$ in the expression
$S^*({\bf t}) = \sum_{i=1}^N|t^{v_i}|$.  Next, observe that
$$\mu(\{{\bf t} \in A_l: S^*({\bf t}) < \epsilon\}) = \int_{\{{\bf t} \in A_l:\,  S^*({\bf t}) < \epsilon\}}\prod_{k=1}^m |{\bf t}_k|^{-\alpha_k} \,d{\bf t} \eqno (2.1)$$
Note  that on $A_l$ a function $|{\bf t}_k|$ is comparable
in magnitude to $u_{b_k}$ for the $b_k$ which is maximal amongst the $u_i$ appearing in ${\bf t}_k$. Thus in terms of the $u$ variables there 
are positive constants $C_1$ and $C_2$ such that  
$$ C_1\int_{\{{\bf u}:\, 0 < u_1 < ... < u_n < r,\,\,S_l^*({\bf u}) < \epsilon\}}\prod_{k=1}^m u_{b_k}^{-\alpha_k} \,d{\bf u} < \mu(\{{\bf t} \in A_l:  S^*({\bf t}) < \epsilon\}) $$
$$< C_2 \int_{\{{\bf u}:\, 0 < u_1 < ... < u_n < r,\,\,S_l^*({\bf u}) < \epsilon\}}\prod_{k=1}^m u_{b_k}^{-\alpha_k} \,d{\bf u} \eqno (2.2)$$
We make a change of variables on $A_l$, setting $u_k = \prod_{i=k}^n y_i$. In the $y$ variables, $A_l$  becomes the rectangular box $(0,1)^{n-1} \times (0,r)$ and on $A_l$,
the function $S_l^*({\bf u})$ becomes $S_l^{**}({\bf y}) = \sum_{i=1}^N y^{W_i}$ for some multiindices $W_i$. Then for some $\beta_k$ the
 integral in $(2.2)$ becomes
$$\int_{\{{\bf y} \in (0,1)^{n-1} \times (0,r):\, S_l^{**}({\bf y}) < \epsilon\}}\prod_{k=1}^n y_k^{-\beta_k} \,d{\bf y} \eqno (2.3)$$ 
Because $K({\bf t})$ is integrable on a neighborhood of the origin, each $\beta_k < 1$ here. We next change variables from $y_k$ to $z_k= y_k^{1 - \beta_k}$
in $(2.3)$. Given $W_i = (W_{i1},...,W_{in})$, let $X_i$ denote $({W_{i1} \over 1 - \beta_1},..., {W_{in} \over 1 - \beta_n})$ and let $S_l^{***}({\bf z}) =
\sum_{i=1}^N z^{X_i}$. Then for some positive constant $C'$ equation  $(2.3)$ becomes
$$C'\int_{\{{\bf z} \in  (0,1)^{n-1} \times (0,r^{1 - \beta_k}):\, S_l^{***}({\bf z}) < \epsilon\}} 1\,d{\bf z} \eqno (2.4)$$
In summary, $\mu(\{{\bf t} \in A_l:  S^*({\bf t}) < \epsilon\})$ is comparable in magnitude to 
the Lebesgue measure of  the set of points in $(0,1)^{n-1} \times (0,r^{1 - \beta_k})$ where $S_l^{***}({\bf z}) < \epsilon$.  We can apply Theorem 1.2 of [G2] to $S_l^{***}({\bf z})$ and say that  for $r > 0$ sufficiently small there exist positive constants $c_{r,l}$, $d_{r,l}$, and $a_l$, and an integer $0 \leq d_l \leq n-1$ such that
$$ c_{r,l}\epsilon^{a_l} |\ln\epsilon|^{d_l} < m(\{{\bf z} \in (0,r)^n: S_l^{***}({\bf z}) < \epsilon\}) < d_{r,l}\epsilon^{a_l} |\ln\epsilon|^{d_l} \eqno (2.5)$$
Here $m$ denotes Lebesgue measure.  By dilation invariance of the existence of such estimates, $(2.5)$ will also hold with $(0,r)^n$ replaced by 
$(0,1)^{n-1} \times (0,r^{1 - \beta_k})$, with different constants.
Technically, Theorem 1.2 of [G2] requires that the components ${W_{ij} \over 1 - \beta_j}$ all be 
integers but the proof of that theorem is valid in the current setting. It should be pointed out that $(2.5)$ can also be shown using related considerations from [V].
 Adding $(2.5)$ over all $l$ gives $(1.5)$, taking $a_0$ to be the minimal $a_l$ and $d_0$ to be the maximal $d_l$ over all those $a_l$ for which $a_l = a_0$.
This concludes the proof of Lemma 2.1.

By Theorem 1.2 of [G1], the quantity $a_l$ is exactly the reciprocal of the Newton distance (Definition 1.3) of $S_l^{***}({\bf z})$, and $d_l$ is $n-1-l$, where $l$ 
denotes the 
 minimal dimension of
any compact face of the Newton polyhedron of $S_l^{***}({\bf z})$ containing the point where the line $z_1 = z_2 = .... = z_n$ intersects this Newton
polyhedron. Here 
a vertex of a Newton polyhedron is taken to be a face of dimension zero. We refer to chapter 7 of [AGuV] for more information.
Thus one can explicitly determine the quantities $a_0$ and $d_0$ using the
Newton polyhedron of $S({\bf x})$ itself, by finding $S^*({\bf x})$, then doing the coordinate changes of the proof of Lemma 2.1 and finding the various Newton 
polyhedra of the functions $S_l^{***}({\bf z})$. 

\noindent We also will make use of the following lemma from [G1].

\noindent {\bf Lemma 2.2. (Lemma 3.7 of [G1]).} Suppose $f({\bf x})$ is a smooth function on a neighborhood of the origin in $\R^n$ with $f(0) = 0$ and $\nabla f(0) = (0,...,0)$. Suppose $f({\bf x})$ has a nonvanishing 
Taylor expansion at the origin, and that either $f({\bf x})$ is real-analytic or a smooth function whose Newton polyhedron intersects
each coordinate axis. 

Then there is a neighborhood $U$ of the origin and constants $K$ and $\eta$ such that if 
$R$ is any dyadic rectangle in $U$ then $R$ may be divided into at most $K$ rectangles $R_j$ such that if
$2^{-k_i}$ denotes the length of $R$ in the $x_i$ direction, for each $R_j$ there is an $a$ and a single
$y =  (y_1,...,y_n)$ with  $|y_i| \leq 2^{-k_i}$ for all $i$ such that on $R_j$ we have
$$|(y \cdot \nabla)^a f({\bf x})| \geq \eta f^*(R_j) \eqno (2.6)$$
Here $f^*(R_j)$ denotes $\sup_{R_j} f^*({\bf x})$. For any $j$, in $(2.6)$ one can arrange that $0 \leq a 
\leq o(f)$ or that $2 \leq a \leq \max(2, o(f))$.

\noindent {\heading 3. Proof of part 1 of Theorem 1.1.} 

\noindent We embed $T$ in an analytic family, defining $T_z$ as follows, where  $S^*({\bf t})$ is as in $(1.4)$.
$$T_z f({\bf x}, x_{n+1}) = e^{z^2}\int_{\R^n} (\min(|t_1|,...,|t_n|, |S^*({\bf t})|))^{z} f({\bf x} - {\bf t}, x_{n+1} - S({\bf t})) K({\bf t})\,d{\bf t}\eqno (3.1)$$
Note that $T_0 = T$. We will show that  if $0 > s_0 > \max(-g, -{1 \over \max(o(S),2)})$, for $z$ on the line $Re(z) = s_0$, if  $2 < p < \infty$ one has estimates $||T_z f||_{L^p}
\leq C_{s_0}||f||_{L^p}$. We will then show that if $s_1 > \max(0,{1 \over \max(o(S),2)} - g)$, then on the line  $Re(z) = s_1$ one has estimates  $||T_z f||_{L_\beta^2}
\leq C_{s_1}||f||_{L^2}$ for $\beta = {1 \over \max(o(S),2)} $. Interpolating between these two estimates gives a Sobolev space estimate for $T = T_0$. Then 
letting $p$ go to infinity, $s_0$ go to $ \max(-g, -{1 \over \max(o(S),2)})$, and $s_1$ go to $\max(0,{1 \over \max(o(S),2)} - g)$, the statement of part 1 of Theorem 1.1 will follow.

\noindent We start by proving the $L^p$ boundedness estimates for $p > 2$.

\noindent {\bf Lemma 3.1.}  Suppose $ 0 > s_0 > \max(-g, -{1 \over \max(o(S),2)})$. Then for $z$ on the line $Re(z) = s_0$, if  $2 < p < \infty$ there exists a constant $C_{s_0}$ such that
 $||T_z f||_{L^p} \leq C_{s_0}||f||_{L^p}$.
 
\noindent {\bf Proof.} Let $s_0 > - g$ and suppose $Re(z) = s_0$. Then $T_z f$ can be expressed as $f \ast \nu_z$ for some measure $\nu_z$, and
as a result for any $2 < p < \infty$,  $||T_z||_{L^p \rightarrow L^p}$ is bounded by the integral of the 
magnitude of the density function in $(3.1)$. So we have the following, where we write $z = s_0 + it$.
$$||T_z||_{L^p \rightarrow L^p} \leq e^{s_0^2 - t^2}\int_{\R^n} (\min(|t_1|,...,|t_n|, |S^*({\bf t})|))^{s_0} |K({\bf t})|\,d{\bf t}\eqno (3.2)$$
 Inserting $(1.2a)$, for some $r_0 > 0$ we have
 $$||T_z||_{L^p \rightarrow L^p} \leq Ce^{s_0^2}\int_{(-r_0,r_0)^n} (\min(|t_1|,...,|t_n|, |S^*({\bf t})|))^{s_0} \prod_{k=1}^m |{\bf t}_k|^{-\alpha_k}\,d{\bf t}\eqno (3.3)$$
 Since $s_0 < 0$, this is in turn bounded by
 $$ Ce^{s_0^2}\int_{(-r_0,r_0)^n} |S^*({\bf t})|^{s_0} \prod_{k=1}^m |{\bf t}_k|^{-\alpha_k}\,d{\bf t} + Ce^{s_0^2}\sum_{i=1}^n \int_{(-r_0,r_0)^n} |t_i|^{s_0} \prod_{k=1}^m |{\bf t}_k|^{-\alpha_k}\,d{\bf t} \eqno (3.4)$$
We first estimate the left term in $(3.4)$. Assuming as we may that $r_0$ is sufficiently small, this term is bounded by
$$Ce^{s_0^2}\sum_{i=0}^{\infty} \int_{\{{\bf t} \in (-r_0,r_0)^n:\,  2^{-i-1} \leq |S^*({\bf t})| < 2^{-i}\}}2^{-is_0}|{\bf t}_k|^{-\alpha_k}\,d{\bf t} \eqno (3.5)$$
Using $(1.5)$ (and recalling that $S^*({\bf t})$ is even in each variable) we see that this is bounded by $C'e^{s_0^2}\sum_{i=0}^{\infty} 2^{-i s_0} \times i^{d_0}2^{-ia_0}$. Since we are assuming that $s_0 > -g \geq -a_0$,
we have that $s_0 + a_0 > 0$ and this sum converges with a bound depending on $s_0$. Thus the first term in $(3.4)$ satisfies the needed bounds.

\noindent Moving on to bounding a given term of the sum on the right of $(3.4)$, let $k_0$ be such that $t_i$ is one of the coordinate functions in ${\bf t}_{k_0}$. Then integrating the term in the the $t$ variables not part of ${\bf t}_{k_0}$, we get a bound of
$$ C \int_{ (-r_0,r_0)^{l_{k_0}}}|t_i|^{s_0} |{\bf t}_{k_0}|^{-\alpha_{k_0}}\,d{\bf t}_{k_0} \eqno (3.6)$$
Performing the integral in the ${\bf t}_{k_0}$ variables other than $t_i$ by breaking  into $|{\bf t}_{k_0}| < 2|t_i| $ and $|{\bf t}_{k_0}| >2|t_i|$ portions, we obtain
$$C'\int_{ (-r_0,r_0)^{l_{k_0}}}|t_i|^{s_0} \max(|t_i|^{l_{k_0} - \alpha_{k_0}-1},1)\,dt_i \eqno (3.7)$$
(If $l_{k_0} - \alpha_{k_0}-1 = 0$ we get an additional $\log|t_i|$ factor that does not affect our argument.)
So we have two cases, depending on whether or not $l_{k_0} - \alpha_{k_0}-1 < 0$. If this does hold, since
$s_0 > -g \geq -(l_{k_0} - \alpha_{k_0})$, the integral in $(3.7)$ converges. On the other hand if $l_{k_0} - \alpha_{k_0}-1 \geq 0$, since $s_0 > -{1 \over \max(o(S),2)} \geq -{1 \over 2}$, the integral once again converges. This concludes the proof of Lemma 3.1.

\noindent {\bf Lemma 3.2.} Suppose $s_1 > \max(0,{1 \over \max(o(S),2)} - g)$. For $z$ on the line $Re(z) = s_1$  there exists a constant $C_{s_1}$ such 
that $||T_z f||_{L_\beta^2}
\leq C_{s_1}||f||_{L^2}$ where $\beta = {1 \over \max(o(S),2)} $.

\noindent {\bf Proof.} Note that $T_z f$ is of the form $f \ast \nu_z$, so to prove Lemma 3.2 one must that $|\hat{\nu}_z(\lambda)| \leq 
C_{s_1}(1 + |\lambda|)^{-\beta}$. Since the integrand in $(3.8)$ below is integrable with an $L^1$ norm bounded by a function of $s_1$, it suffices to prove the estimate
$|\hat{\nu}_z(\lambda)| \leq C_{s_1}|\lambda|^{-\beta}$ for say $|\lambda| > 2$.

\noindent Explicitly, we have 
$$\hat{\nu}_z(\lambda) = e^{z^2}\int_{\R^n} e^{i\lambda_1 t_1 + i\lambda_2 t_2 +... +i\lambda_n t_n + i\lambda_{n+1}S({\bf t})}(\min(|t_1|,...,|t_n|, |S^*({\bf t})|))^{z} K({\bf t})\,d{\bf t}\eqno (3.8)$$
We write $\hat{\nu_z}(\lambda) = \sum_j I_j$, where for $j = (j_1,...,j_n)$, $I_j$ denotes the portion of the integral $(3.8)$ over the dyadic rectangle $R_j =
\{{\bf t} \in \R^n: 2^{-j_i - 1} \leq |t_i| < 2^{-j_i}$ for each $i\}$. By Lemma 2.2, we  can further write $I_j$ as the sum of boundedly many integrals over 
rectangles on which $(2.6)$ holds for the function $S({\bf t})$. Each such rectangle can in turn be written as the sum of boundedly many subsets on which the minimum in 
$\min(|t_1|,...,|t_n|, |S^*({\bf t})|)$ is achieved by a specific $|t_1|,...,|t_n|,$ or $|S^*({\bf t})|$. We label the subsets corresponding to a given $I_j$ as $R_{jk}$ and we correspondingly write $(3.8)$ as the sum of
$I_{jk}$, where 
$$I_{jk} = e^{z^2}\int_{R_{jk}} e^{i\lambda_1 t_1 + i\lambda_2 t_2 +... + i\lambda_n t_n +i\lambda_{n+1}S({\bf t})}(\min(|t_1|,...,|t_n|, |S^*({\bf t})|))^{z} K({\bf t})\,d{\bf t}\eqno (3.9)$$
Regardless of what $\lambda$ is, we may select a $\lambda_i$ such that $|\lambda_i| > {1 \over n+1}|\lambda|$. We will have two separate arguments,
depending on whether $i < n + 1$ or $i = n+1$. We first suppose  that $i < n + 1$.  

\noindent {\bf Case 1.} $|\lambda_i| > {1 \over n+1}|\lambda|$ for some $i \leq n$.

Let 
$P({\bf t})$ denote the phase function in $(3.9)$. Note that $\partial_{t_i} P({\bf t}) = \lambda_i + \lambda_{n+1} \partial_{t_i} S({\bf t})$. Since 
$S({\bf t})$ has a zero of order at least two at the origin by $(1.3)$, we can assume that we are in a sufficiently small neighborhood of the origin such that
$$|\partial_{t_i} P({\bf t})| > {1 \over 2n + 2} |\lambda| \eqno (3.10)$$ 
We focus our attention on the $t_i$ integration in $(3.9)$ for the other variables fixed. There are boundedly many intervals of integration, and on each 
we integrate by parts, integrating $e^{\lambda_1 t_1 + \lambda_2 t_2 +... + \lambda_n t_n +\lambda_{n+1}S({\bf t})}$
and differentiating the rest. Due to $(3.10)$, we gain a factor of $C{1 \over |\lambda|}$ in the integration by parts, but due to $(1.2b)$ and the 
corresponding derivative estimates on $|S^*({\bf t})|^z$ we also get a factor of $C|z| {1 \over |t_i|}$ from the differentiation.
As a result, substituting in $(1.2a)$ we have
$$|I_{jk}| \leq C|ze^{z^2}|  \int_{R_j}{1 \over |\lambda t_i|}|(\min(|t_1|,...,|t_n|, |S^*({\bf t})|))^{z}| \prod_{k=1}^m |{\bf t}_k|^{-\alpha_k}\,d{\bf t}\eqno (3.11a)$$
(The endpoint terms in the integration by parts will satisfy the same bounds.)
Since $|ze^{z^2}|$ is bounded on any vertical line in the complex plane, the above is bounded by
$$= C_{s_1}  \int_{R_j}{1 \over |\lambda t_i|}(\min(|t_1|,...,|t_n|, |S^*({\bf t})|))^{s_1} \prod_{k=1}^m |{\bf t}_k|^{-\alpha_k}\,d{\bf t}\eqno (3.11b)$$
Using that $s_1 > 0$, we may take $|t_i|$ in the minimum and therefore $(3.11b)$ is bounded by 
$$ C_{s_1}'\int_{R_j}{1 \over |\lambda t_i|}|t_i|^{s_1} \prod_{k=1}^m |{\bf t}_k|^{-\alpha_k}\,d{\bf t}\eqno (3.12)$$
By simply taking absolute values and integrating in $(3.9)$, we also have a bound of
$$|I_{jk}| \leq C_{s_1}''\int_{R_j}|t_i|^{s_1} \prod_{k=1}^m |{\bf t}_k|^{-\alpha_k}\,d{\bf t}\eqno (3.13)$$
Combining $(3.12)$ and $(3.13)$ gives
$$|I_{jk}| \leq C_{s_1}'''\int_{R_j}\min\bigg(1, {1 \over |\lambda t_i|}\bigg)|t_i|^{s_1} \prod_{k=1}^m |{\bf t}_k|^{-\alpha_k}\,d{\bf t}\eqno (3.14)$$
Adding $(3.14)$ over all $j$ and $k$ gives for some $r_0 > 0$ that 
$$|\hat{\nu}_z(\lambda)| \leq C_{s_1}'''\int_{(-r_0,r_0)^n}\min\bigg(1, {1 \over |\lambda t_i|}\bigg)|t_i|^{s_1} \prod_{k=1}^m |{\bf t}_k|^{-\alpha_k}\,d{\bf t}\eqno (3.15)$$
Integrating $(3.15)$ in the $t$ variables other than $t_i$, similarly to $(3.7)$ we get
$$|\hat{\nu}_z(\lambda)| \leq C_{s_1}''''\int_0^{r_0}\min\bigg(1, {1 \over |\lambda| t_i}\bigg)t_i^{s_1} \max(t_i^{l_i - \alpha_i - 1}, 1)\,dt_i\eqno (3.16)$$
Once again if $l_i - \alpha_i - 1$ happens to equal zero we get an additional $|\ln t_i|$ that does not affect our arguments. It is natural to break up
the integral of $(3.16)$ into $t_i < {1 \over |\lambda|}$ and $t_i > {1 \over |\lambda|}$ portions. The first part is given by
$$\int_0^{{1 \over |\lambda|}}t_i^{s_1} \max(t_i^{l_i - \alpha_i - 1}, 1)\,dt_i\eqno (3.17)$$
Since $s_1 > 0$, if one is taking $1$ in the maximum one obtains a bound of $C|\lambda|^{-1}$ which is far better than we need. So assume that 
$l_i - \alpha_i - 1 < 0$ and we are taking $t^{l_i - \alpha_i - 1}$ in the maximum. 
Since this lemma assumes that $s_1 > {1 \over \max(o(S),2)} - g$ and by definition of $g$ we have $g \leq l_i - \alpha_i$, we have that 
$s_1 + l_i - \alpha_i > {1 \over \max(o(S),2)} - g + l_i - \alpha_i > {1 \over \max(o(S),2)}$. Then integrating $(3.17)$ gives a bound of $C|\lambda|^{-(s_1 + l_i - \alpha_i)} \leq C|\lambda|^{-{1 \over \max(o(S),2)} - \epsilon_{s_1}}$ for some $\epsilon_{s_1} > 0$, giving the bounds needed in this lemma.

Moving on to the $t_i > {1 \over |\lambda|}$ portion of $(3.16)$, note that the magnitude of the part of the integral over the
dyadic pieces $t_i \sim 2^{-j}$ will increase or decrease
exponentially in $j$, with the integral over the outermost dyadic piece being of order $|\lambda|^{-1}$, which is better than the estimate needed. Given
that the integral over the innermost dyadic piece has the same bound of $C|\lambda|^{-{1 \over \max(o(S),2)} - \epsilon_{s_1}}$ as above,
the sum of these magnitudes will again have a bound of  $C|\lambda|^{-{1 \over \max(o(S),2)} - \epsilon_{s_1}}$. This completes the argument for 
the case where $|\lambda_i| > {1 \over n+1}|\lambda|$ for some $i \leq n$.

\noindent {\bf Case 2.} $|\lambda_{n+1}| > {1 \over n+1}|\lambda|$. 

\noindent We will make use of the Van der Corput lemma (see p. 334 of [S]).

\noindent {\bf Van der Corput lemma.} Suppose $P(x)$ is a real-valued $C^k$ function on the interval $[a,b]$ with $|P^{(k)}(x)| > M$ on $[a,b]$ for
some $M > 0$. Let $\psi(x)$ be a complex-valued $C^1$ function on $[a,b]$. If $k \geq 2$ there is a constant $c_k$ depending only on $k$ such that
$$\bigg|\int_a^b e^{iP(x)}\psi(x)\,dx\bigg| \leq c_kM^{-{1 \over k}}\bigg(|\psi(b)| + \int_a^b |\psi'(x)|\,dx\bigg) \eqno (3.18)$$
If $k =1$, the same is true if we add the conditions that $P(x)$ is $C^2$ and that $P'(x)$ is monotonic on $[a,b]$. 

We examine the function $S({\bf t})$ on the set $R_{jk}$ of $(3.9)$. Recall that $R_{jk}$ is contained in the dyadic rectangle $R_j =
 \{{\bf t} \in \R^n: 2^{-j_i - 1} \leq |t_i| < 2^{-j_i}$ for each $i\}$. Let $S_1({\bf t}) = S(2^{-j_1}t_1,...,2^{-j_n}t_n)$ and $S_1^*({\bf t})=  S^*(2^{-j_1}t_1,...,2^{-j_n}t_n)$.
Then by Lemma 2.2 there is a constant $\eta > 0$, a direction $v$ and an 
$a$ satisfying $2 \leq a \leq \max(o(S),2)$ such that on $R_{jk}^* = \{(t_1,...,t_n): (2^{-j_1}t_1,...,2^{-j_n}t_n) \in R_{jk}\}$ for any fixed
${\bf t}_0$ in $R_j$ one has the estimate
$$|\partial_v^a S_1({\bf t})| > \eta |S_1^*({\bf t}_0)| \eqno (3.19)$$
Note that in the dilated coordinates, if $K_1({\bf t})$ denotes $K(2^{-j_1}t_1,...,2^{-j_n}t_n)$, then equation $(3.9)$ becomes
$$I_{jk} = 2^{-(j_1 + ... + j_n)} e^{z^2}\int_{R_{jk}^*} e^{i2^{-j_1}\lambda_1 t_1 + i2^{-j_2}\lambda_2 t_2 +... + i2^{-j_n}\lambda_n t_n + i\lambda_{n+1}S_1({\bf t})}$$
$$\times  (\min(2^{-j_1}|t_1|,...,2^{-j_n}|t_n|, |S_1^*({\bf t})|))^{z}K_1({\bf t})\,d{\bf t}\eqno (3.20)$$
We apply the Van der Corput lemma for $a$th derivatives  in the $v$ direction in $(3.20)$, and then integrate the result in the $n-1$ orthogonal directions. Note that 
due to the way the $R_j$ were divided into $R_{jk}$, there may be boundedly many intervals of integration in the $v$ direction, in which case we apply the
Van der Corput lemma on each interval and add the result. 

In the situation at hand the quantity denoted by $M^{-{1 \over k}}$ in $(3.18)$ is given by $|\lambda_{n+1} S_1^*({\bf t}_0)|^{-{1 \over a}}$. Since $|\lambda_{n+1}| > {1 \over n +1}|\lambda|$ in our argument, this in turn is bounded by $C|\lambda S_1^*({\bf t}_0)|^{-{1 \over  a}}$. In view of the 
$|\psi'(x)|$ in the integral $(3.18)$, we examine the effect of a derivative in the $v$ direction on the various factors in $(3.20)$. Since $R_{jk} \subset [{1 \over 2},1]^n$, a derivative landing on $(\min(2^{-j_1}|t_1|,...,2^{-j_n}|t_n|, |S_1^*({\bf t})|))^{z}$ will introduce a factor of $C|z|$. In view of
$(1.2a)-(1.2b)$ and the fact that $R_{jk} \subset [{1 \over 2},1]^n$, a derivative landing on $K_1({\bf t})$ introduces a factor of $C$ when using $(1.2a)-
(1.2b)$ in upper bounds. Thus applying the Van der Corput lemma in the $v$ direction has the overall effect of introducing a factor bounded by 
$C|z||\lambda S_1^*({\bf t}_0)|^{-{1 \over  a_0}}$. Translating this back into the unscaled coordinates, we get
$$I_{jk} \leq C|z e^{z^2}| \int_{R_j} |\lambda S^*({\bf t}_0)|^{-{1 \over a}}\min(|t_1|,...,|t_n|, |S^*({\bf t})|)^{s_1}\prod_{k=1}^m |{\bf t}_k|^{-\alpha_k}\,d{\bf t} \eqno (3.21)$$
Note that $|z e^{z^2}|$ is uniformly bounded on a given line $Re\, z = s_1$, and that $|S^*({\bf t})|$ varies by a factor of at most $C$ on any dyadic
rectangle. Thus $(3.21)$ implies
$$I_{jk} \leq C_{s_1}\int_{R_j} |\lambda S^*({\bf t})|^{-{1 \over a}}\min(|t_1|,...,|t_n|, |S^*({\bf t})|)|^{s_1}\prod_{k=1}^m |{\bf t}_k|^{-\alpha_k}\,d{\bf t} \eqno (3.22a)$$
Adding this over the boundedly many $k$ for a given $j$, we get the same form for the estimate for $I_j$, the portion of the integral $(3.8)$ coming from the 
rectangle $R_j$:
$$I_j \leq C_{s_1}\int_{R_j} |\lambda S^*({\bf t})|^{-{1 \over a}}\min(|t_1|,...,|t_n|, |S^*({\bf t})|)^{s_1}\prod_{k=1}^m |{\bf t}_k|^{-\alpha_k}\,d{\bf t} \eqno (3.22b)$$
By simply taking absolute values of the integrand and integrating, we also have
$$I_j \leq C_{s_1}'\int_{R_j} \min(|t_1|,...,|t_n|, |S^*({\bf t})|)^{s_1}\prod_{k=1}^m |{\bf t}_k|^{-\alpha_k}\,d{\bf t} \eqno (3.23)$$
Combining $(3.22b)$ and $(3.23)$ and using $|S^*({\bf t})|$ in the minimum (note that we are using that $s_1 \geq 0$ here), we get
$$I_j \leq C_{s_1}''\int_{R_j} \min(1, |\lambda S^*({\bf t})|^{-{1 \over a}})|S^*({\bf t})|^{s_1}\prod_{k=1}^m |{\bf t}_k|^{-\alpha_k}\,d{\bf t} \eqno (3.24)$$
Since $a \leq \max(o(S),2)$, this in turn is bounded by
$$I_j \leq C_{s_1}''\int_{R_j} \min(1, |\lambda S^*({\bf t})|^{-{1 \over \max(o(S),2)}})|S^*({\bf t})|^{s_1}\prod_{k=1}^m |{\bf t}_k|^{-\alpha_k}\,d{\bf t} \eqno (3.25)$$
Finally, adding $(3.25)$ over all dyadic rectangles $R_j$ gives for some $r_0 > 0$ that
$$|\hat{\nu}_z(\lambda)| \leq C_{s_1}'''\int_{(-r_0,r_0)^n} \min(1, |\lambda S^*({\bf t})|^{-{1 \over \max(o(S),2)}})|S^*({\bf t})|^{s_1}\prod_{k=1}^m |{\bf t}_k|^{-\alpha_k}\,d{\bf t} \eqno (3.26a)$$
Since everything in the integral $(3.26a)$ is even in all variables, we can replace $(-r_0,r_0)^n$ by $(0,r_0)^n$ and say
$$|\hat{\nu}_z(\lambda)| \leq C_{s_1}''''\int_{(0,r_0)^n} \min(1, |\lambda S^*({\bf t})|^{-{1 \over \max(o(S),2)}})|S^*({\bf t})|^{s_1}\prod_{k=1}^m |{\bf t}_k|^{-\alpha_k}\,d{\bf t} \eqno (3.26b)$$
In terms of the measure $d\mu = \prod_{k=1}^m |{\bf t}_k|^{-\alpha_k}\,d{\bf t}$, this can be rewritten as
$$|\hat{\nu}_z(\lambda)| \leq C_{s_1}''''\int_{(0,r_0)^n} \min(1, |\lambda S^*({\bf t})|^{-{1 \over \max(o(S),2)}})|S^*({\bf t})|^{s_1}\,d\mu \eqno (3.26c)$$
It is natural to break up $(3.26c)$ into $|S^*({\bf t})| < |\lambda|^{-1}$ and $|S^*({\bf t})| > |\lambda|^{-1}$ portions. We obtain
$$|\hat{\nu}_z(\lambda)| \leq C_{s_1}'''' \int_{|S^*({\bf t})| < |\lambda|^{-1}}|S^*({\bf t})|^{s_1}\,d\mu$$
$$ + \,\,\, C_{s_1}''''
\int_{|S^*({\bf t})| > |\lambda|^{-1}}|\lambda|^{-{1 \over \max(o(S),2)}}|S^*({\bf t})|^{s_1 - {1 \over \max(o(S),2)} }\,d\mu\eqno (3.27)$$
If one writes the first integral in $(3.27)$ as the sum in $i$ of integrals over the points where $2^{-i-1} |\lambda|^{-1} <  |S^*({\bf t})| \leq 2^{-i} |\lambda|^{-1}$, and inserts $(1.5)$ into each term, since $s_1 > 0$ one gets a geometric sum that decreases as $i$ increases. Thus the overall sum can be bounded by a
constant times the $i = 0$ term, which by $(1.5)$ is bounded by $C|\lambda|^{-s_1 - a_0}(\ln |\lambda|)^{d_0}$. Since we are assuming $s_1 > 
{1 \over  \max(o(S),2)} - g$ and since $g \leq a_0$, we have $s_1 + a_0 > {1 \over  \max(o(S),2)}$ and as a result the bound 
$C|\lambda|^{-s_1 - a_0}(\ln |\lambda|)^{d_0}$ implies a bound of $C'|\lambda|^{-{1 \over  \max(o(S),2)} -\epsilon_{s_1}}$ for some $\epsilon_{s_1} > 0$.

Similarly, one can write  the second integral as the sum in $i$ of integrals over the points ${\bf t}$ where $2^{i} |\lambda|^{-1} <  |S^*({\bf t})| \leq 2^{i + 1} |\lambda|^{-1}$ and then insert $(1.5)$. The resulting estimates change exponentially in $i$, with the $i = 0$ term being comparable to $C|\lambda|^{-s_1 - a_0}(\ln |\lambda|)^{d_0} < C'|\lambda|^{-{1 \over  \max(o(S),2)} -\epsilon_{s_1}}$  like above, while the last term will be comparable to $C|\lambda|^{-{1 \over \max(o(S),2)}}$
As a result we have a bound of $C''|\lambda|^{-{1 \over \max(o(S),2)}}$ for the whole second integral in $(3.27)$. Combining with the above bounds for
the first integral, we see that we have an estimate of the form $|\hat{\nu}_z(\lambda)| \leq C_{s_1}'''''|\lambda|^{-{1 \over \max(o(S),2)}}$, as is 
needed for Lemma 3.2. This concludes the argument for the case where $|\lambda_{n+1}| > {1 \over n + 1}|\lambda|$ and therefore
the proof of Lemma 3.2.

\noindent {\bf The end of the proof of part 1 of Theorem 1.1.}

By Lemma 3.2, for $Re\,(z) = s_1 > \max({1 \over  \max(o(S),2)} - g,0)$ we have $||T_z f||_{L^2_{\beta}} \leq C_{s_1}||f||_{L^2}$ with $\beta = 
{1 \over \max(o(S),2)}$. Furthermore, by Lemma 3.1, if $2 < p < \infty$ and $Re\,(z) = s_0 > \max(-g, -{1 \over \max(o(S),2)})$, then we have an estimate
$||T_z f||_{L^p} \leq C_{s_0}||f||_{L^p}$.

 Note that
$0 = \alpha\max(-g, -{1 \over \max(o(S),2)}) +  (1 - \alpha)(\max({1 \over  \max(o(S),2)} - g,0)) $, where $\alpha = \max(0, 1 - g\max(2,o(S)))$. 
Thus  if $1 > \alpha' > \alpha$,
 one can write $0 = \alpha' s_0 + (1 - \alpha') s_1$, where $0 > s_0 >\max(-g, -{1 \over \max(o(S),2)})$ and  $s_1  > \max({1 \over  \max(o(S),2)} - g,0)$. Hence by  complex interpolation $T = M_0$ is bounded from $L^q$ to $L^q_{\beta}$, where
 ${1 \over q} =  \alpha' {1 \over p} + (1 -  \alpha') {1 \over 2}$ and $\beta = \alpha' 0 + (1 - \alpha'){1 \over \max(o(S),2)}$. Explicitly, we have
 $q = {1 \over {1 \over 2} + \alpha'({1 \over p} - {1 \over 2})}$ and $\beta = {1 - \alpha' \over \max(o(S),2)}$. 
 
 Using interpolation again, we have that $T$ is bounded from $L^r$ to $L^r_{\gamma}$ for $({1 \over r}, \gamma)$ in the closed 
 triangle with vertices $(0,0), (1,0),$ and $({1 \over q}, \beta)$. Taking the union of these triangles as $\alpha'$ approaches $\alpha$ and $p$ approaches
 $\infty$, we get that
 $T$ is bounded from $L^r$ to $L^r_{\gamma}$  for $({1 \over r}, \gamma)$ in the open triangle with vertices $(0,0), (1,0),$ and $({1 \over q'}, \beta')$, where 
$q' =  {1 \over {1 \over 2} - {1 \over 2}\alpha} = {2 \over \min(1,g\max(o(S), 2))} = \max(2, {2 \over g\max(o(S),2)})$ and where $\beta' = 
 {1 - \alpha \over  \max(o(S),2)} = \min({1 \over  \max(o(S), 2)}, g)$.
 
 In the case where $g \geq {1 \over  \max(o(S), 2)}$, the union of these triangles is the open triangle with vertices $(0,0)$, $(1,0)$, and $({1 \over 2},
 {1 \over  \max(o(S), 2)})$, which is the boundedness region stipulated by part 1 of Theorem 1.1 in this case.
 If $g < {1 \over  \max(o(S), 2)}$, 
 the union of these triangles  is the open triangle with vertices $(0,0)$, $(1,0)$, and $({g \max(o(S),2) \over 2}, g)$. So $T$ is bounded from $L^r$ to $L^r_{\gamma}$ for $({1 \over r}, \gamma)$ in this region. By duality, it is also  bounded from $L^r$ to $L^r_{\gamma}$ for $({1 \over r}, \gamma)$ 
 such that $(1 - {1 \over r}, \gamma)$ is in this region, giving the triangle with vertices $(0,0)$, $(1,0)$, and $(1 - {g \max(o(S),2) \over 2}, g)$. 
 Thus $T$ is bounded from $L^r$ to $L^r_{\gamma}$ for $({1 \over r}, \gamma)$ in the open trapezoidal region with vertices $(0,0)$, $(1,0)$, 
 $({g \max(o(S),2) \over 2}, g)$, and $(1 - {g \max(o(S),2) \over 2}, g)$. This is the region stipulated by Theorem 1.1 in the case where  $g < {1 \over  \max(o(S), 2)}$. This completes the proof of part 1 of Theorem 1.1.

\noindent  {\heading 4. Proof of part 2 of Theorem 1.1.}

 Assume the hypotheses of part 2 of Theorem 1.1 hold. Namely, assume that $g < 1$, $K({\bf t})$ is nonnegative, and there exists a positive constant $C_0$ and a neighborhood $N_0$ of the origin such that $K({\bf t}) > C_0\prod_{k=1}^m |{\bf t}_k|^{-\alpha_k}$ on $\{{\bf t} \in N_0: t_i \neq 0$ for all $i\}$. Suppose that  $T$ is bounded from $L^p(\R^{n+1})$ to $L^p_{\beta}(\R^{n+1})$ for some $ 1 < p < \infty$. Then by duality, $T$ is bounded from $L^q(\R^{n+1})$ to $L^q_{\beta}(\R^{n+1})$ where $q$ is such that ${1 \over p} + {1 \over q} = 1$. Since either $p \leq 2 \leq q$ or $q \leq 2 \leq p$, interpolation shows that
$T$ is bounded from $L^2(\R^{n+1})$ to $L^2_{\beta}(\R^{n+1})$. As a result, if $\nu$ is the measure such that $Tf = f \ast \nu_z$, we have an estimate
$|\hat{\nu}(\lambda)| \leq C(1 + |\lambda|)^{-\beta}$. Explicitly, we have
$$\bigg|\int_{\R^n} e^{i\lambda_1 t_1 + i\lambda_2 t_2 +... + i\lambda_n t_n + i\lambda_{n+1}S({\bf t})} K({\bf t})\,d{\bf t}\bigg| \leq C(1 + |\lambda|)^{-\beta}
\eqno (4.1)$$
We will first show that $\beta \leq a_0$, and then
we will show that $\beta \leq l_i - \alpha_i$ for each $i$. Since $g = \min(a_0, l_1 - \alpha_1,...,\l_m - \alpha_m)$, this will give part 2 of Theorem 1.1.
Since $(4.1)$ holds in all directions, it holds in the $(0,...,0,\lambda_{n+1})$ direction, so $(4.1)$ implies that for all $\lambda_{n+1}$ we have
$$\bigg|\int_{\R^n} e^{i\lambda_{n+1}S({\bf t})} K({\bf t})\,d{\bf t}\bigg| \leq C(1 + |\lambda_{n+1}|)^{-\beta} \eqno (4.2)$$
Denote the integral on the left of $(4.2)$ by $U(\lambda_{n+1})$.
Let $B(x)$ be a bump function on $\R$  whose Fourier transform is nonnegative, compactly supported, and equal to 1 on a 
neighborhood of the origin, and let $\epsilon$ be a small  positive number. If $0 <  \beta' < \beta$, then $(4.2)$ implies
 that for some constant $A$ independent of $\epsilon$ one has
$$\int_{\R}|U(\lambda_{n+1}) \lambda_{n+1}^{\beta '  - 1}B(\epsilon \lambda_{n+1})|\,d \lambda_{n+1} < A \eqno (4.3)$$
As a result we have
$$\bigg|\int_{\R^{n+1}} e^{i\lambda_{n+1}S({\bf t})} K({\bf t}) |\lambda_{n+1}|^{\beta'  - 1}B(\epsilon \lambda_{n+1})\,d\lambda_{n+1} \,d{\bf t}\bigg|  <  A \eqno (4.4)$$
We do the integral in $\lambda_{n+1}$ in $(4.4)$. Letting $b_N(y)$ be the convolution of $|y|^{-\beta '}$ with ${1 \over \epsilon} \hat{B}({y \over \epsilon})$,  for a constant $A'$ independent of $\epsilon$ we get
$$\bigg|\int_{\R^n}b_{\epsilon} (-S({\bf t})) \,K({\bf t})\,d{\bf t}\bigg| < A' \eqno (4.5)$$
Note that both $b_{\epsilon} (-S({\bf t}))$ and $K({\bf t})$ are nonnegative here. Thus we may remove the absolute value and let
 $\epsilon \rightarrow 0$ to obtain
$$\int_{\R^n}|S({\bf t})|^{-\beta '} K({\bf t}) < \infty  \eqno (4.6)$$
Since $ K({\bf t})$ is bounded below by $C_0\prod_{k=1}^m |{\bf t}_k|^{-\alpha_k}$ on a neighborhood $N_0$ of the origin, we therefore have
$$\int_{N_0}|S({\bf t})|^{-\beta '}\prod_{k=1}^m |{\bf t}_k|^{-\alpha_k} \,d{\bf t}< \infty  \eqno (4.7)$$
In other words, $|S({\bf t})|^{-\beta '}$ is in $L^1(N_0)$ with respect to the measure $\mu$. Hence it is in weak $L^1$, 
and we have the existence of a constant $C$ such that
$$\mu (\{{\bf t} \in N_0 : |S({\bf t})|^{-\beta '} > \epsilon \}) \leq C  {1 \over  \epsilon} \eqno (4.8)$$
Replacing $ \epsilon$ by $ \epsilon^{-\beta '}$, we get
$$\mu (\{{\bf t} \in N_0: |S({\bf t})| <  \epsilon \}) \leq C  \epsilon^{\beta '} \eqno (4.9)$$
By Lemma 2.1 of [G2] there is a constant $C'$ such that $|S({\bf t})| \leq C'S^*({\bf t})$, so we also have
$$\mu (\{{\bf t} \in N_0: S^*({\bf t}) <  \epsilon\}) \leq C''\epsilon^{\beta '} \eqno (4.10)$$
In view of the definition of $a_0$, we have $\beta ' \leq a_0$. Since this holds for each $\beta'$ satisfying
$0 < \beta' < \beta$, we conclude that $\beta \leq a_0$ as needed.

Showing that $\beta \leq l_i - \alpha_i$ for each $i$ is quite similar. This time we use $(4.1)$ in the $(0,...,0,\lambda_i,0,...,0)$ direction, and the steps from
$(4.2)$ to $(4.7)$ lead to
$$\int_{N_0}|t_i|^{-\beta '}\prod_{k=1}^m |{\bf t}_k|^{-\alpha_k} \,d{\bf t}< \infty  \eqno (4.11)$$
Suppose $k_0$ is such that $t_i$ is a component of ${\bf t}_{k_0}$. Then integrating $(4.11)$ in the remaining variables first leads to the following holding for
some $\delta > 0$.
$$\int_{(0,\delta)^{l_i}} |t_i|^{-\beta '} |{\bf t}_{k_0}|^{-\alpha_{k_0}} \,d{\bf t}_{k_0}< \infty  \eqno (4.12)$$
As a result, $\beta' + \alpha_{k_0} < l_{k_0}$. Since this is true for all $0 < \beta' < \beta$, we conclude that $\beta \leq l_{k_0} - \alpha_{k_0}$.
Since $t_i$ was arbitrary, this holds for all $k_0$. This completes the proof of part 2 of Theorem 1.1.

\noindent {\heading 5. References.}

\noindent [AGuV] V. Arnold, S. Gusein-Zade, A. Varchenko, {\it Singularities of differentiable maps},
Volume II, Birkhauser, Basel, 1988. \parskip = 3pt\baselineskip = 3pt

\noindent [C] M. Christ, {\it Failure of an endpoint estimate for integrals along curves} in Fourier analysis and partial differential equations (Miraflores de la Sierra, 1992), 163-168, Stud. Adv. Math., CRC, Boca Raton, FL, 1995. 

\noindent [CNSW] M. Christ, A. Nagel, E. M. Stein, and S. Wainger, {\it Singular and maximal Radon transforms: analysis and geometry},
Ann. of Math. (2) {\bf 150} (1999), no. 2, 489-577. 

\noindent [Cu] S. Cuccagna, {\it  Sobolev estimates for fractional and singular Radon transforms}, J. Funct. Anal. {\bf 139} (1996), no. 1, 94-118.

\noindent [Gr] L. Grafakos, {\it Endpoint bounds for an analytic family of Hilbert transforms}, Duke Math. J. {\bf 62} (1991), no. 1, 23-59. 

\noindent [G1] M. Greenblatt, {\it Maximal averages over hypersurfaces and the Newton polyhedron}, J. Funct. Anal. {\bf 262} (2012), no. 5, 2314-2348. 

\noindent [G2] M. Greenblatt, {\it Oscillatory integral decay, sublevel set growth, and the Newton
polyhedron}, Math. Annalen {\bf 346} (2010), no. 4, 857-895.

\noindent [G3] M. Greenblatt, {\it Smooth and singular maximal averages over 2D hypersurfaces  and associated Radon transforms}, submitted

\noindent [G4] M. Greenblatt, {\it Uniform bounds for Fourier transforms of surface measures in $\R^3$ with nonsmooth density},
Trans. Amer. Math. Soc. {\bf 368} (2016), no. 9, 6601-6625.

\noindent [G5] M. Greenblatt, {\it An analogue to a theorem of Fefferman and Phong for averaging operators
along curves with singular fractional integral kernel}, Geom. Funct. Anal. {\bf 17} (2007), no. 4, 1106-1138.

\noindent [IKM] I. Ikromov, M. Kempe, and D. M\"uller, {\it Estimates for maximal functions associated
to hypersurfaces in $\R^3$ and related problems of harmonic analysis}, Acta Math. {\bf 204} (2010), no. 2,
151--271.

\noindent [Ka1] V. N. Karpushkin, {\it A theorem concerning uniform estimates of oscillatory integrals when
the phase is a function of two variables}, J. Soviet Math. {\bf 35} (1986), 2809-2826.

\noindent [Ka2] V. N. Karpushkin, {\it Uniform estimates of oscillatory integrals with parabolic or 
hyperbolic phases}, J. Soviet Math. {\bf 33} (1986), 1159-1188.

\noindent [PSe] M. Pramanik, A. Seeger, {\it $L^p$ Sobolev regularity of a restricted X-ray transform in $\R^3$} (English summary) Harmonic analysis and its applications, 47-64, Yokohama Publ., Yokohama, 2006. 

\noindent [Se] A. Seeger, {\it Radon transforms and finite type conditions}, J. Amer. Math. Soc. {\bf 11} (1998), no. 4, 869-897.

\noindent [S] E. Stein, {\it Harmonic analysis; real-variable methods, orthogonality, and oscillatory
integrals}, Princeton Mathematics Series Vol. 43, Princeton University Press, Princeton, NJ, 1993.

\noindent [St1] B. Street, {\it Sobolev spaces associated to singular and fractional Radon transforms}, Rev. Mat. Iberoam. {\bf 33} (2017), no. 2, 633-748.

\noindent [St2] B. Street, {\it Multi-parameter Singular Integrals}, Annals ofMathematics Studies, {\bf 189}, Princeton
University Press, Princeton, NJ, 2014.

\noindent [V] A. N. Varchenko, {\it Newton polyhedra and estimates of oscillatory integrals}, Functional 
Anal. Appl. {\bf 18} (1976), no. 3, 175-196.
\end